\documentstyle[12pt,epsf]{article}
\def\mysection#1{\refstepcounter{section}\subsection{#1}}
\def\mysubsection#1{\subsubsection{#1}}

\def\ns{\normalsize}
\def\ssize{\normalsize}
\def\note#1{}

\def\ds{\displaystyle}
\def\ts{\textstyle}

\newcommand{\nn}{\nonumber}

\def\al{\alpha}
\def\rh{\rho}
\def\del{\delta}
\def\la{\lambda}
\def\m{\mu}
\def\ta{\tau}
\def\l{\left}
\def\r{\right}
\def\oti{\otimes}
\def\un{{1\kern-.25em{\rm l}}}
\def\tri{\triangle}

\begin{document}

\begin{titlepage}
\rightline{NIKHEF 95-059}
\vskip 1.8  true cm
\begin{center}
\Large{\bf Link invariants from $N$-state vertex models: an
alternative construction independent of statistical models}
\\
\vspace{0.45in}
\ns\sc
M.J. Rodr\'\i guez-Plaza
\\
\vspace{0.3in}
\ssize\em
NIKHEF,
Postbus 41882,
1009 DB Amsterdam, The Netherlands\\
\vspace{0.04in}
and\\
\vspace{0.04in}
\ns\sc
\ssize\em
Institut f\"ur Theoretische Physik,
Universit\"at Heidelberg,
D-69120 Heidelberg, Germany
\footnote{current address}\\
\end{center}
\vspace{0.8in}
{\leftskip=1.5 true cm \rightskip=1.5 true cm
\noindent
We reproduce the hierarchy of link invariants associated to the series
of $N$-state vertex models with a method different from the original
construction due to Akutsu, Deguchi and Wadati. The alternative method
substitutes the `crossing symmetry' property exhibited by the
Boltzmann weights of the vertex models by a similar property which,
for the purpose of constructing link invariants, encodes the same
information but requires only the limit of the Boltzmann
weights when the spectral parameter is sent to infinity.  \par}
\vskip .8 true cm
\end{titlepage}
\setcounter{page}{2}

\mysection{Introduction}
\label{intro}

Starting from the $N$-state vertex models first introduced in
\cite{SAA}, Akutsu, Deguchi and Wadati show in \cite{ADW} that there
is a polynomial link invariant associated to each vertex model of the
series. The invariant corresponds to a Markov trace and is
therefore a link invariant of ambient isotopy for oriented links.
In particular for $N=2$ (the 6-vertex model) the {\em skein relation}
of the polynomial link invariant is given by
\begin{equation}
\al\l({{}\atop\epsfbox{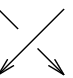}}\r)=
(1-t)\,t^{1/2}\,\al\l({{}\atop\epsfbox{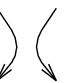}}\r)
+t^2\,\al\l({{}\atop\epsfbox{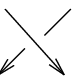}}\r)
\label{polN2}
\end{equation}
that corresponds to the Jones's polynomial \cite{Jo}. For $N=3$ (the
19-vertex model) it is given by
\begin{eqnarray}
\al\l({{}\atop\epsfbox{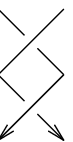}}\r)&=&
t\,(1-t^2+t^3)\,\al\l({{}\atop\epsfbox{sigma1.eps}}\r)\nn\\
&&+t^2\,(t^2-t^3+t^5)\,\al\l({{}\atop\epsfbox{sigma0.eps}}\r)-
t^8\,\al\l({{}\atop\epsfbox{sigma-1.eps}}\r),
\label{polN3}
\end{eqnarray}
which is a one-variable specialization of the Kauffman polynomial
\cite{Kau1}; for $N=4$ (the 44-vertex model) the
relation is
\begin{eqnarray}
\al\l({{}\atop\epsfbox{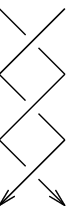}}\r)&=&
t^{3/2}\,(1-t^3+t^5-t^6)\,\al\l({{}\atop\epsfbox{sigma2.eps}}\r)\nn\\
&&+ t^6\,(1-t^2+t^3+t^5-t^6+t^8)\,\al\l({{}\atop\epsfbox{sigma1.eps}}\r)
\nn\\
&&- t^{9/2}\,t^8\,(1-t+t^3-t^6)\,\al\l({{}\atop\epsfbox{sigma0.eps}}\r)-
t^{20}\,\al\l({{}\atop\epsfbox{sigma-1.eps}}\r),
\label{polN4}
\end{eqnarray}
that is again a one-variable polynomial.  This sequence generalizes
for $N$ arbitrary in a polynomial defined by a $N$-th order skein
relation.

The object of this paper is to prove that the same $N=2, 3, 4$
polynomials (\ref{polN2})-(\ref{polN4}) are obtained with the link
invariant that we recall next. The case of generic $N$ can be
worked out by induction and its associated link invariant again
reproduces the result obtained in \cite{ADW}.  The link invariant is
the following. Consider the plane projection of any classical link so
that the projected link diagram consists only of double crossings,
maxima, minima and vertical arcs and associate to each of these pieces
the next objects indexed by a finite index set $I$
\[
{{}\atop\epsfbox{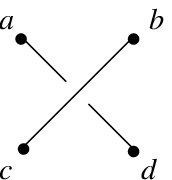}}\longleftrightarrow\,\,R^a{}_c{}^b{}_d,
\qquad\qquad{{}\atop\epsfbox{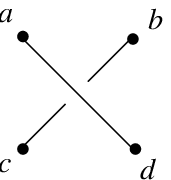}}\longleftrightarrow
\,\,{R^{-1}}^a{}_c{}^b{}_d,
\]
\[{{}\atop\epsfbox{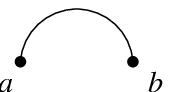}}\longleftrightarrow\,\,M_{a\,b},
\qquad\qquad{{}\atop\epsfbox{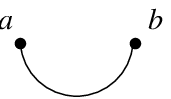}}\longleftrightarrow\,\,M^{a\,b},
\qquad\qquad
{{}\atop\epsfbox{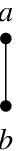}}\longleftrightarrow\,\,\delta^a{}_b.
\]
With this convention any link diagram $L$ is translated into its
corresponding expression $<L>$ in terms of the previous elements. Thus
we have, for example, that for the trefoil 
\[
{{}\atop\epsfbox{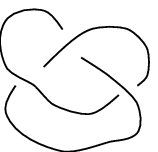}}\qquad\longleftrightarrow\qquad 
{{}\atop\epsfbox{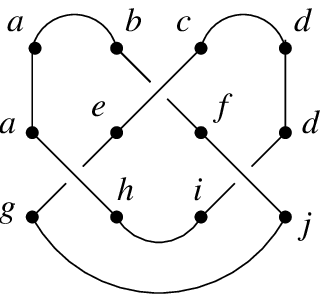}}
\]
\[ <\, {\rm trefoil}\,>= M_{a\,b}\,M_{c\,d}\,R^b{}_e{}^c{}_f\,
{R^{-1}}^a{}_g{}^e{}_h\,{R^{-1}}^f{}_i{}^d{}_j\,M^{h\,i}\,M^{g\,j}
\]
and for the unknotted circle
${{}\atop\epsfbox{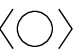}}=M_{a\,b}\,M^{a\,b}$, where sum over
repeated indices is always assumed.

It is well-known \cite{Kau2} \cite{KR} that if the objects
$R^a{}_c{}^b{}_d$, ${R^{-1}}^a{}_c{}^b{}_d$, $M_{a\,b}$, $M^{a\,b}$
have been chosen so that they satisfy the conditions
\begin{eqnarray}
&&M_{a\,b}\,M^{b\,c}=\del_a{}^c=M^{c\,b}\, M_{b\,a},\qquad\quad
{{}\atop\epsfbox{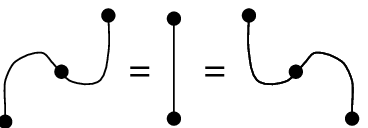}}\label{m}\\
&&R^a{}_c{}^b{}_d\, {R^{-1}}^c{}_e{}^d{}_f=\del_e{}^a\,\del_f{}^b=
{R^{-1}}^a{}_c{}^b{}_d\,R^c{}_e{}^d{}_f,\qquad\quad
{{}\atop\epsfbox{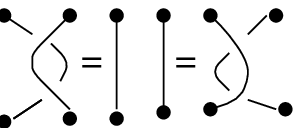}}\label{r}\\
&&R^a{}_i{}^b{}_j\, R^j{}_k{}^c{}_f\, R^i{}_d{}^k{}_e=
R^b{}_i{}^c{}_j\, R^a{}_d{}^i{}_k\, R^k{}_e{}^j{}_f,
\qquad\quad{{}\atop\epsfbox{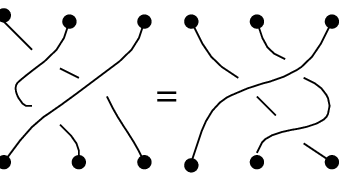}}\label{braid}\\
&&{R^{-1}}^a{}_c{}^b{}_d=M^{a\,e}\,R^b{}_e{}^f{}_c\,M_{f\,d},
\qquad\quad{{}\atop\epsfbox{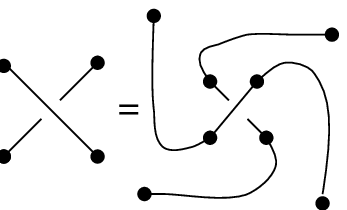}}\label{twist1}\\
&&{R^{-1}}^a{}_c{}^b{}_d=M_{c\,e}\,R^e{}_d{}^a{}_f\,M^{f\,b},
\qquad\quad{{}\atop\epsfbox{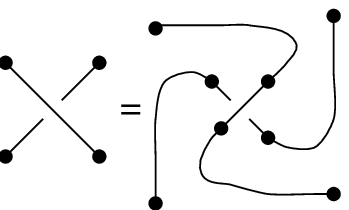}}\label{twist2}\\
\end{eqnarray}
then $<L>$ is an invariant of regular isotopy for unoriented link
diagrams.  Equations (\ref{m}) and (\ref{r}) require that the matrices
$M_u=(M^{a\,b})$ and $M_d=(M_{a\,b})$  are inverse to each other, and
also that $R$ and $R^{-1}$ are (actually the notation used anticipated
it already); (\ref{braid}) is the {\em Yang-Baxter relation} to which
$R$ has to be a solution (we will refer as $R$-matrix to any solution of
this equation). The mixed conditions (\ref{twist1}) and (\ref{twist2})
are `crossing symmetry'-like properties demanded for $R$.

The paper is organized as follows. First we solve eqs
(\ref{twist1})-(\ref{twist2}) using as initial input the $R$-matrices
written in (\ref{R2}), (\ref{R3}) and (\ref{R4}) respectively, where
$Z$ is a constant to be determined. The solution to these equations
gives the matrices $M_u$ and $M_d$ corresponding to each input $R$ and
constitutes the main result of Section~\ref{two}.  Section~\ref{three}
lists sufficient conditions for the link invariant $<L>$ to behave as
a Markov trace. In this and the following section $L$, the link
diagram, is represented as the closure of a braid since braids are
algebraically more manageable than projections of links.
Section~\ref{four} checks that for the three set of matrices $R$,
$M_u$ and $M_d$ under consideration, the corresponding link invariant
$<L>$ behaves as a Markov trace and constructs the ambient isotopy
invariant that derives from such Markov trace. The result are
precisely the link invariants (\ref{polN2}), (\ref{polN3}) and
(\ref{polN4}) obtained in \cite{ADW} through the $N=2,3,4$ state
vertex model. Our conclusions are redacted in
Section~\ref{conclusions}. There are other methods known to provide
(part or the whole series of) $N$-state invariants such as Kauffman
bracket approach and state models \cite {Kau1} \cite{Kau3} and
Kirillov and Reshetikhin work \cite{KR}. Section~\ref{six} explains
the relation between $<L>$ and these methods.  The last section
contains some remarks that we find of interest.

The way in which Akutsu {\em et al} derived their result is different
from the method that we follow here. The construction of a Markov
trace in \cite{ADW} relies on the crossing-symmetry property exhibited
by the matrix $R(u)$ of Boltzmann weights of the solvable $N$-state
vertex model, and on the non-trivial expression of the crossing
multipliers. To deduce these, the crossing multipliers, reference
\cite{ADW} uses the explicit dependence of $R(u)$ on the spectral
parameter $u$. In the method followed here is not essential to know
how $R(u)$ depends on $u$ because the information contained in the
crossing multipliers useful to write a Markov trace is substituted by
the information encoded in eqs (\ref{twist1}) and (\ref{twist2}). It
was precisely the similarity between the latter two equations with the
crossing symmetry property (\ref{crossing}) what originally motivated
this work.  To reproduce the link polynomials
(\ref{polN2})-(\ref{polN4}) then, there exist other methods in
addition to vertex models (via $<L>$ for instance) which reproduce
independently the results originally obtained with statistical models.

The origin of matrices written in (\ref{R2}), (\ref{R3}) and
(\ref{R4}) is this: they are the limit $R=Z\,\lim_{u\to\infty}
R(u)/\rh(u)$ taken on the matrix $R(u)$ constructed, as we have said,
with the Boltzmann weights of the $N=2,3,4$ vertex model. The
denominator $\rh(u)$ is a function evaluated with each particular
$R(u)$. At the expense of being repetitive we emphasize once more that
in this paper we are working not with $R(u)$ but with its limit and
that the limit is enough to write the invariant $<L>$. 

\mysection{Solving the $N=2,3,4$ cases}
\label{two}
\mysubsection{$N=2$ case}
Let us consider the matrix $R(u)$
constructed with the Boltzmann weights of the $N=2$ vertex model
as in \cite{ADW}
\begin{eqnarray}
R(u)=\l(\begin{array}{cccc}
\sinh(\la-u)&0&0&0\\
0&e^{\ds{2\,\m\,u}}\sinh\la&\sinh u&0\\
0&\sinh u&e^{\ds{-2\,\m\,u}}\sinh\la&0\\
0&0&0&\sinh(\la-u)
\label{spectN2}
\end{array}\r)
\end{eqnarray}
(along this paper we are using $R^i{}_l{}^j{}_k(u)$ to denote what in
\cite{ADW} is denoted by $S^i{}_j{}^k{}_l(u)$).  Here $u$ represents
the spectral parameter of the vertex model and $\la$ and $\m$ are
arbitrary constants. This matrix is a solution of the Yang-Baxter
equation with spectral parameter
\begin{eqnarray}
R^a{}_i{}^b{}_j(u)\, R^j{}_k{}^c{}_f(u+v)\, R^i{}_d{}^k{}_e(v)=
R^b{}_i{}^c{}_j(v)\, R^a{}_d{}^i{}_k(u+v)\, R^k{}_e{}^j{}_f(u)
\label{YBE}
\end{eqnarray}
where $a,b,c,\ldots$ are indices in the set $I=\{-1/2,\,1/2\}$, and
it satisfies the relation $R(u)\,R(-u)=\rh(u)\,\rh(-u)$, with
$\rh(u)\equiv\sinh(\la-u)$. It also satisfies the crossing symmetry
property
\begin{equation}
R^i{}_k{}^j{}_l(u)=\l({r(i)\,r(k)\over r(j)\,r(l)}\r)^{1/2}
R^j{}_{-i}{}^{-l}{}_k(\la-u)
\label{crossing}
\end{equation}
with crossing multipliers given by
$r(p)=e^{\ds{-2\,\m\,\la\,p}}$ for every $p$ in $I$.
The limit $R=Z\,\lim_{u\to\infty} R(u)/\rh(u)$ is well defined and
given by the invertible matrix with expression
\begin{eqnarray}
R=Z\pmatrix{
1&0&0&0\cr
0&1-q^2&q&0\cr
0&q&0&0\cr
0&0&0&1\cr}
\label{R2}
\end{eqnarray}
in terms of the parameter $q\equiv -e^\la$. This matrix here is
obtained with the assumption $\m=1/2$ in (\ref{spectN2}) and satisfies
the Yang-Baxter relation (\ref{braid}).

With the $R$-matrix (\ref{R2}) and its inverse we now determine
the value of the constant $Z$ with the condition that there exist
non-singular matrices $M_u$ and $M_d$, inverse to each other, such
that they satisfy conditions (\ref{twist1}) and (\ref{twist2}).  It is
not difficult to see that for generic $q$ these equations have a
unique solution given by
\begin{eqnarray}
M_u=\l(\begin{array}{cc}
0     &    q^{1/2}      \\
-q^{-1/2}  &    0
\end{array}\r),
\qquad M_d=- M_u,\qquad  Z=\pm\, q^{-1/2},
\label{M2}
\end{eqnarray}
where we have used the freedom to choose a multiplicative constant in
$M_u$, say, to make both matrices of determinant equal to one. This is
done for simplicity merely because the value of the multiplicative
constant is not relevant since it does not affect the link invariant
$<L>$ derived from the solutions of (\ref{twist1})-(\ref{twist2}).
This solution yields for the zero knot the value of the invariant
${{}\atop\epsfbox{zero.eps}}=M_{a\,b}\,M^{a\,b}={\rm
tr}\,(M_u\,M_d^t)=-(q+q^{-1})$ where the superscript $t$ indicates
transpose matrix. We postpone until Section~\ref{four} the
final expression of the link invariant associated to (\ref{R2}) and
(\ref{M2}).

\mysubsection{$N=3$ case} 
In this case the $R$-matrix considered as input 
to solve eqs (\ref{twist1}) and (\ref{twist2}) is
given by the invertible matrix
\begin{eqnarray}
R=Z\l(\begin{array}{ccccccccc}
1&0&0&0&0&0&0&0&0\\
0&1-q^4&0&-q^2&0&0&0&0&0\\
0&0&(1-q^2)\,(1-q^4)&0&q\,(1-q^4)&0&q^4&0&0\\
0&-q^2&0&0&0&0&0&0&0\\
0&0&q\,(1-q^4)&0&q^2&0&0&0&0\\
0&0&0&0&0&1-q^4&0&-q^2&0\\
0&0&q^4&0&0&0&0&0&0\\
0&0&0&0&0&-q^2&0&0&0\\
0&0&0&0&0&0&0&0&1
\end{array}\r).
\,\label{R3}
\end{eqnarray}
This matrix is $Z$ times the $\lim_{u\to\infty} R(u)/\rh(u)$ where now
$R(u)$ is the spectral parameter dependent solution of (\ref{YBE})
with entries the Boltzmann weights of the $N=3$ vertex model given
by \cite{ADW}
\begin{eqnarray}
&&R^1{}_1{}^1{}_1(u)=R^{-1}{}_{-1}{}^{-1}{}_{-1}(u)=
\sinh(\la-u)\,\sinh(2\,\la-u),\nn\\
&&R^1{}_{-1}{}^{-1}{}_1(u)=R^{-1}{}_1{}^1{}_{-1}(u)=
\sinh u\,\sinh(\la+u),\nn\\
&&e^{\ds{4\,\m\,u}}\,R^1{}_1{}^{-1}{}_{-1}(u)=
e^{\ds{-4\,\m\,u}}\,R^{-1}{}_{-1}{}^1{}_1(u)=
\sinh \la\,\sinh 2\,\la,\nn\\
&&R^1{}_0{}^0{}_1(u)=R^{-1}{}_0{}^0{}_{-1}(u)=
R^0{}_1{}^1{}_0(u)=R^0{}_{-1}{}^{-1}{}_0(u)=
\sinh u\,\sinh(\la-u),\nn\\
&&e^{\ds{2\,\m\,u}}\,R^1{}_1{}^0{}_0(u)=
e^{\ds{-2\,\m\,u}}\,R^{-1}{}_{-1}{}^0{}_0(u)=
e^{\ds{-2\,\m\,u}}\,R^0{}_0{}^1{}_1(u)
\label{N3}\\
&&\qquad\qquad\qquad =e^{\ds{2\,\m\,u}}\,R^0{}_0{}^{-1}{}_{-1}(u)=
\sinh 2\,\la\,\sinh(\la-u),\nn\\
&&e^{\ds{-2\,\m\,u}}\,R^0{}_{-1}{}^0{}_1(u)=
e^{\ds{2\,\m\,u}}\,R^0{}_1{}^0{}_{-1}(u)=
e^{\ds{2\,\m\,u}}\,R^1{}_0{}^{-1}{}_0(u)\nn\\
&&\qquad\qquad\qquad = e^{\ds{-2\,\m\,u}}\,R^{-1}{}_0{}^1{}_0(u)
=\sinh 2\,\la\,\sinh u,\nn\\
&&R^0{}_0{}^0{}_0(u)=\sinh \la\,\sinh 2\,\la-\sinh u\,\sinh(\la-u)\nn
\end{eqnarray}
and zero the rest of the entries. In this case the index set is fixed
as $I=\{-1,\,0,\, 1\}$; the quantities $R^a{}_c{}^b{}_d$ are arranged
in matrix form so that the left indices label the block and the right
ones the block entries. Solution (\ref{N3}) satisfies
$R(u)\,R(-u)=\rh(u)\,\rh(-u)$ where now
$\rh(u)\equiv\sinh(\la-u)\,\sinh(2\,\la-u)$ is different from the
$N=2$ case but obtained with the same relation, and also satisfies
property (\ref{crossing}) from above with the same crossing
multipliers. As in the case $N=2$, the matrix in (\ref{R3}) is
obtained from (\ref{N3}) with the choice $\mu=1/2$ and the
substitution $q=-e^\la$.

For generic values of $q$ there is a unique solution (again unique up
to a multiplicative constant) to eqs
(\ref{twist1}) and (\ref{twist2}) corresponding to the matrix
in (\ref{R3}) and its inverse. The solution is given by
\begin{eqnarray}
M_u=\l(\begin{array}{ccc}
0  & 0   &   q    \\
0  & -1  &    0   \\
q^{-1} & 0 & 0
\end{array}\r),\qquad M_d=M_u, \qquad  Z=\pm\, q^{-2}.
\label{M3}
\end{eqnarray}
The invariant associated to the unknot is in this case
equal to ${\rm tr}\,(M_u\,M_d^t)=q^2+1+q^{-2}$.

\mysubsection{$N=4$ case} 
The non-zero entries of the $R$-matrix 
$R=Z\,\lim_{u\to\infty} R(u)/\rh(u)$ where $R(u)$ \cite{ADW} are the
Boltzmann weights of the $N=4$ vertex model are given by (now
$\rh(u)\equiv\sinh(\la-u)\,\sinh(2\,\la-u)\,\sinh(3\,\la-u)$ and
$\mu=1/2$, $q=-e^\la$ as usual)
\begin{eqnarray}
&&R^{3/2}{}_{3/2}{}^{3/2}{}_{3/2}=
R^{-3/2}{}_{-3/2}{}^{-3/2}{}_{-3/2}=Z\nn\\
&&R^{-3/2}{}_{3/2}{}^{3/2}{}_{-3/2}=
R^{3/2}{}_{-3/2}{}^{-3/2}{}_{3/2}=Z\,q^9\nn\\
&&R^{-3/2}{}_{-3/2}{}^{3/2}{}_{3/2}=Z\,(1-q^2)\,(1-q^4)\,(1-q^6)\nn\\
&&R^{3/2}{}_{1/2}{}^{1/2}{}_{3/2}=
R^{-3/2}{}_{-1/2}{}^{-1/2}{}_{-3/2}=
R^{1/2}{}_{3/2}{}^{3/2}{}_{1/2}=
R^{-1/2}{}_{-3/2}{}^{-3/2}{}_{-1/2}=Z\,q^3\nn\\
&&R^{1/2}{}_{-3/2}{}^{-3/2}{}_{1/2}=
R^{-1/2}{}_{3/2}{}^{3/2}{}_{-1/2}=
R^{-3/2}{}_{1/2}{}^{1/2}{}_{-3/2}=
R^{3/2}{}_{-1/2}{}^{-1/2}{}_{3/2}=Z\,q^6\nn\\
&&R^{-3/2}{}_{-3/2}{}^{-1/2}{}_{-1/2}=
R^{1/2}{}_{1/2}{}^{3/2}{}_{3/2}=Z\,(1-q^6)\nn\\
&&R^{1/2}{}_{-3/2}{}^{-1/2}{}_{3/2}=
R^{-3/2}{}_{1/2}{}^{3/2}{}_{-1/2}=Z\,q^4\,(1-q^6)\label{R4}\\
&&R^{-3/2}{}_{-1/2}{}^{1/2}{}_{-1/2}=
R^{1/2}{}_{-1/2}{}^{1/2}{}_{3/2}=
R^{-1/2}{}_{1/2}{}^{3/2}{}_{1/2}=
R^{-1/2}{}_{-3/2}{}^{-1/2}{}_{1/2}\nn\\
&&\qquad\qquad\qquad=Z\,q^3\,(1-q^4)\,(q^2+1+q^{-2})^{1/2}\nn\\
&&R^{-3/2}{}_{-3/2}{}^{1/2}{}_{1/2}=
R^{-1/2}{}_{-1/2}{}^{3/2}{}_{3/2}=Z\,(1-q^4)\,(1-q^6)\nn\\
&&R^{-1/2}{}_{-3/2}{}^{1/2}{}_{3/2}=
R^{-3/2}{}_{-1/2}{}^{3/2}{}_{1/2}=Z\,q\,(1-q^4)\,(1-q^6)\nn\\
&&R^{1/2}{}_{1/2}{}^{1/2}{}_{1/2}=
R^{-1/2}{}_{-1/2}{}^{-1/2}{}_{-1/2}=Z\,q^4\nn\\
&&R^{1/2}{}_{-1/2}{}^{-1/2}{}_{1/2}=
R^{-1/2}{}_{1/2}{}^{1/2}{}_{-1/2}=Z\,q^5\nn\\
&&R^{-1/2}{}_{-1/2}{}^{1/2}{}_{1/2}=Z\,q^2\,(1-q^4)\,(1+q^2),\nn
\end{eqnarray}
where now $I=\{-3/2,\,-1/2,\,1/2,\,3/2\}$. Using these matrix elements
as input data in (\ref{twist1}) and (\ref{twist2}) the unique solution
(up to a multiplicative constant) to these equations for generic $q$
is
\begin{eqnarray}
M_u=\l(\begin{array}{cccc}
0  & 0   & 0  &  q^{3/2}  \\
0  & 0   & -q^{1/2} &  0   \\
0  &q^{-1/2} & 0 & 0     \\
-q^{-3/2} & 0  &  0 & 0
\end{array}\r),\qquad M_d=-M_u, \qquad  Z=\pm\, q^{-9/2}
\label{M4}
\end{eqnarray}
which gives for the unknot the invariant ${{}\atop\epsfbox{zero.eps}}=
{\rm tr}\,(M_u\,M_d^t)=-(q^3+q+q^{-1}+q^{-3})$.

\mysection{Conditions for the link invariant $<L>$ to be a Markov trace}
\label{three}

We study in this section under which conditions $<L>$ is an invariant
of ambient isotopy in addition to regular isotopy as well. Just for
convenience we will regard every link diagram $L$ as the closure of a
certain braid $A$ in the $n$-string braid group $B_n$. {}From the
picture
\[
{{}\atop\epsfbox{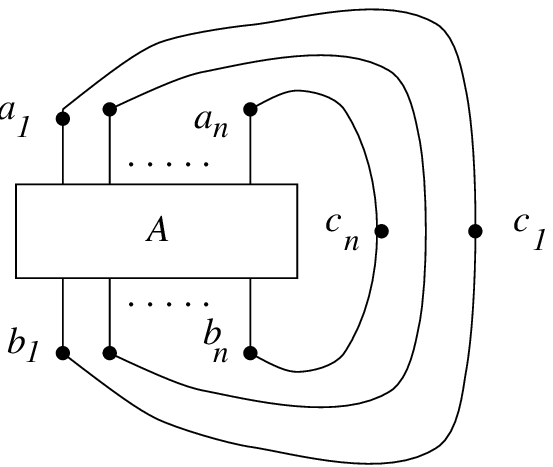}}
\] 
we easily see that the invariant associated to $L$ is  
\begin{eqnarray*}
<L>&=&A^{a_1}{}_{b_1}^{\ldots}{}_{\ldots}^{\,a_n}{}_{b_n}
\,(M^{{b_1}\,{c_1}}\,M_{{a_1}\,{c_1}})
\cdots (M^{{b_n}\,{c_n}}\,M_{{a_n}\,{c_n}})\\
&=&A^{a_1}{}_{b_1}^{\ldots}{}_{\ldots}^{\,a_n}{}_{b_n}
\,(M_u\,M_d^t)^{b_1}{}_{a_1}\cdots (M_u\,M_d^t)^{b_n}{}_{a_n}\\
&=&{\rm tr}\, (A\,(M_u\,M_d^t)^{\otimes^n}).
\end{eqnarray*}
As mentioned $A$ represents an element of the $n$-string braid
group $B_n$ generated by $\{1,\, b_i\}$, $i=1,\ldots, n-1$. The matrix
representation of the elementary braids
\[
{{}\atop\epsfbox{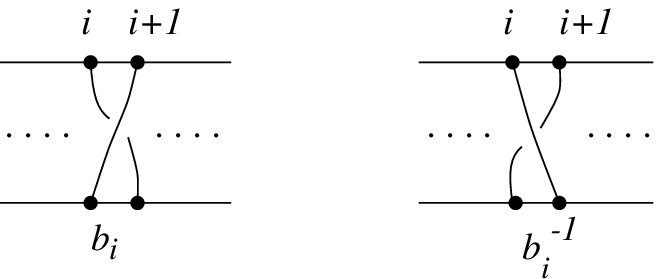}}
\]
that we shall be using to write each $A$ is given in terms of
$R$-matrices by
\begin{equation}
b_i=\underbrace{
\un\oti\cdots\oti R\oti\cdots\un}_{n\; {\rm times}}\,,\qquad\quad
b_i^{-1}=\un\oti\cdots\oti R^{-1}\oti\cdots\un
\label{bi}
\end{equation}
where $R$ or its inverse $R^{-1}$ are placed in the $(i,i+1)$ entry of
the tensor product and are given by (\ref{R2}), (\ref{R3}) or
(\ref{R4}), respectively. In this formula $\un$ denotes the unit
matrix of dimensions $N\times N$. The trace `tr' in $<L>={\rm tr}\,
(A\,(M_u\,M_d^t)^{\otimes^n})$ is the ordinary trace of matrices taken
in this representation. The fact that $R$ satisfies the Yang-Baxter
relation (\ref{braid}) is equivalent to say that the generators
$\{b_i\}$ satisfy $b_i\,b_{i+1}\,b_i=b_{i+1}\,b_i\,b_{i+1}$, $1\le
i\le n-2$ that together with $b_i\,b_j=b_j\,b_i$, $|i-j|\ge 2$ coming
from (\ref{bi}) guarantee the topological equivalence between
different expressions of a braid in terms of the braid generators.

We now define on $B_n$ the functional $\phi:\,B_n\longrightarrow \,
{\cal C}$ by
\begin{equation}
\phi(A)\equiv {{\rm tr}\, (A\,(M_u\,M_d^t)^{\otimes^n})\over
({\rm tr}\, (M_u\,M_d^t))^n}
\label{markov}
\end{equation}
so that $\phi$ is basically $<L>$ but satisfies $\phi(1)=1$.  We
investigate now under which conditions this functional $\phi$ written
as in (\ref{markov}) behaves as a {\em Markov trace}.  For $\phi$ to
be called a Markov trace it must satisfy the following properties
\begin{description}
\item{(p1)} $\phi(A\,B)=\phi(B\,A)$ for all $A,\, B$ in $B_n$ and
\item{(p2)} $\phi(A\,b_n)=\ta\,\phi(A)$ and
$\phi(A\,b_n^{-1})=\bar\ta\,\phi(A)$ for $A$ in $B_n$ and
$b_n,\,b_n^{-1}$ in $B_{n+1}$ with $\ta$ and $\bar\ta$ constants
independent of $n$ and given by
\[
\ta=\phi(b_i),\qquad\bar\ta=\phi(b_i^{-1})\qquad {\rm for}\quad
{\rm all}\quad i.
\]
\end{description}
If $\phi$ is a Markov trace then it is possible to associate an
invariant of ambient isotopy for oriented links $\al'$ to it 
given by the renormalization formula
\begin{equation}
\al'(A)=\l({1\over \ta\,\bar\ta}\r)^{(n-1)/2}\,
\l({\bar\ta\over\ta}\r)^{e(A)/2}\,\phi(A).
\label{ambientp}
\end{equation}
In this formula $e(A)$ is the writhe of the braid $A$ which is given
by the exponent sum of the generators $\{bi\}$ in the braid so that
if, for instance, $A=b_2^3\, b_1^{-2}$ then $e(A)=3-2=1$.

It is not difficult to prove that properties (p1)-(p2) above hold for
$\phi$ given by (\ref{markov}) when the objects $R,\,M_u$ and
$M_d$, subject already to restrictions (\ref{m})-(\ref{twist2}),
satisfy as well the following conditions
\begin{description}
\item{(c1)} \begin{equation}
R\,(M_u\,M_d^t)^{\otimes^2}=(M_u\,M_d^t)^{\otimes^2}\,R,
\label{c1}
\end{equation}
\end{description}
which in its pictorial form means that the crossings can be
pulled through the closure strands and
\begin{description}
\item{(c2)} \begin{equation}
(R\,(M_u\,M_d^t)^{\otimes^2})^a{}_b{}^c{}_c
\cdot\,{\rm tr}\,(M_u\,M_d^t)=
(M_u\,M_d^t)^a{}_b\cdot\,{\rm
tr}\,(R\,(M_u\,M_d^t)^{\otimes^2}),\label{c2}
\end{equation}
\end{description}
together with the similar relations that result from replacing $R$
with $R^{-1}$ in (\ref{c1}) and (\ref{c2}). As above $a,b,c$ are
elements of the index set $I$.
\begin{figure}
\[{{}\atop\epsfbox{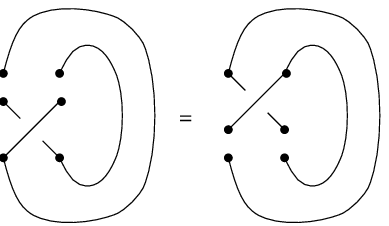}}, \qquad\qquad {{}\atop\epsfbox{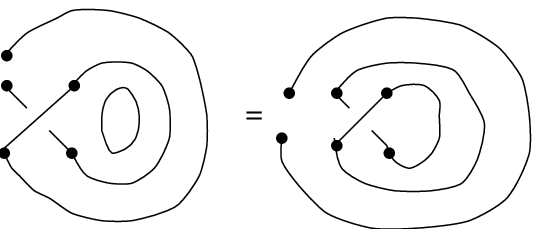}}\]
\caption{Conditions c1 and c2} 
\end{figure}
The proof of conditions (c1), (c2) is as follows. 
For ${\rm tr}\,(A\,B\,(M_u\,M_d^t)^{\otimes^n})= {\rm
tr}\,(B\,A\,(M_u\,M_d^t)^{\otimes^n})$ to hold is enough to demand
that $A\,(M_u\,M_d^t)^{\otimes^n}=(M_u\,M_d^t)^{\otimes^n}\,A$ for any
$A$, which turns to be equivalent to condition (\ref{c1}) since
any braid $A$ decomposes in a product of generators of $B_n$ and these
are expressed in terms of $R$. In the case of condition (\ref{c2})
it merely comes from $\phi$ written for the link
\[
{\epsfbox{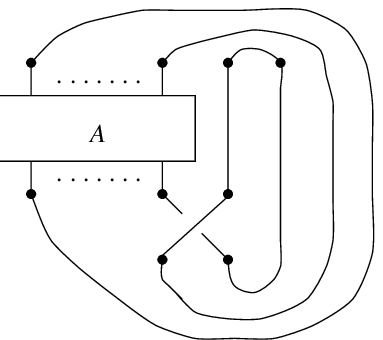}}
\]
and made equal to the product $\phi(A)\,\phi(b_i)$ where
$\phi(b_i)={\rm tr}\,(R\,(M_u\,M_d^t)^{\otimes^2})/({\rm tr}\,
(M_u\,M_d^t))^2$ according to definition (\ref{markov}).

\mysection{Ambient isotopy $N=2,3,4\ldots$ link invariant}
\label{four}

In this section we check that the link invariant $<L>$ behaves as a
Markov trace for the $N=2,3,4$ matrices $R$, $M_u$ and $M_d$, and
calculate the ambient isotopy invariant $\al'(\cdot)$ associated. We
compare the result with the ambient isotopy invariant $\al(\cdot)$
obtained in \cite{ADW} using vertex models.

\mysubsection{$N=2$ case} 
It is simple to check that the matrices $R$, $M_u,\,M_d$ given in
(\ref{R2}) and (\ref{M2}) satisfy conditions (\ref{c1}) and (\ref{c2})
what indicates that $\phi$ defined as in (\ref{markov}) is indeed a
Markov trace. This allows to construct the invariant (\ref{ambientp})
as follows. With the defining formulae
\[
\phi(b_i)={{\rm tr}\,(R\,(M_u\,M_d^t)^{\otimes^2})\over ({\rm tr}\,
(M_u\,M_d^t))^2},\qquad 
\phi(b_i^{-1})={{\rm tr}\,(R^{-1}\,(M_u\,M_d^t)^{\otimes^2})\over ({\rm tr}\,
(M_u\,M_d^t))^2}
\]
and (\ref{R2}), (\ref{M2}) we obtain that 
\[
\ta=\pm\, {q^{-1/2}\over q^2+1},\qquad
\bar\ta=\pm\, {q^{5/2}\over q^2+1}
\]
which substituted in $\al'(\cdot)$ gives
\[
\al'(A)=(q+q^{-1})^{n-1}\,q^{3\,e(A)/2}\,\phi(A).
\]
A consequence of the minimal polynomial of $R$ in (\ref{R2})
\[
(R-Z)\,(R+q^2\,Z)=0,\qquad Z=\pm\, q^{-1/2}
\]
and the linearity of the trace function is the existence of a relation
between the numbers $\phi(b_i),\,\phi(1)$ and $\phi(b_i^{-1})$ that,
due to (\ref{ambientp}), implies a relation between
$\al'(b_i),\,\al'(1)$ and $\al'(b_i^{-1})$ given by
\[
\al'(b_i)=\mp\, q\,(q^2-1)\,\al'(1)+q^{4}\,\al'(b_i^{-1}).
\]
With the substitution $t=q^2$ this equality transforms in
\begin{equation}
\al'\l({{}\atop\epsfbox{sigma1.eps}}\r)=\pm\, 
(1-t)\,t^{1/2}\,\al'\l({{}\atop\epsfbox{sigma0.eps}}\r)
+t^2\,\al'\l({{}\atop\epsfbox{sigma-1.eps}}\r)
\label{polpN2}
\end{equation}
which reproduces the result obtained in \cite{ADW} and displayed in
(\ref{polN2}). The two link invariants calculated with this formula
(note that there is a $(\pm)$ in it) are in fact the same one because
when they are computed for a given knot or link either they do
coincide or differ in a global sign. It can be said then that 
(\ref{polpN2}) defines a unique link invariant. This 
is the Jones polynomial obviously.

\mysubsection{$N=3$ case}
In this case $\phi$ given by (\ref{markov}) is a Markov trace
too since $R$, $M_u$ and $M_d$ in (\ref{R3}) and (\ref{M3}) 
also satisfy (\ref{c1}) and (\ref{c2}). It only remains to proceed as in
the case $N=2$ and compute $\ta$ and $\bar\ta$ which result
\[
\ta=\pm\, {q^{-2}\over q^4+q^2+1},\qquad
\bar\ta=\pm\, {q^{6}\over q^4+q^2+1}
\]
thus providing the invariant
\begin{equation}
\al'(A)=(q^2+1+q^{-2})^{n-1}\,q^{4\,e(A)}\,\phi(A).
\label{ambient3}
\end{equation}
The minimal polynomial of $R$ in (\ref{R3}) given now by
\[
(R-Z)\,(R+q^4\,Z)\,(R-q^6\,Z)=0,\qquad  Z=\pm\, q^{-2}
\]
fixes a relation between $\phi(b_i^2),\,\phi(b_i),\,\phi(1)$ and
$\phi(b_i^{-1})$ that translates into a relation among the $\al'$ on
the same arguments. After some elementary algebra this relation
reads as ($t=q^2$)
\begin{eqnarray}
\al'\l({{}\atop\epsfbox{sigma2.eps}}\r)&=&
\pm\,t\,(1-t^2+t^3)\,\al'\l({{}\atop\epsfbox{sigma1.eps}}\r)\nn\\
\label{polpN3}\\
&&+t^2\,(t^2-t^3+t^5)\,\al'\l({{}\atop\epsfbox{sigma0.eps}}\r)
\mp\,t^8\,\al'\l({{}\atop\epsfbox{sigma-1.eps}}\r),\nn
\end{eqnarray}
that also coincides with formula (\ref{polN3}) as obtained in Akutsu {\em
et al} work. Here again the two ambient link invariant in
(\ref{polpN3}) correspond to the same invariant.

\mysubsection{$N=4$ case} 
For $N=4$ and its associated matrices $R$, $M_u$ and $M_d$ given by
(\ref{R4}) and (\ref{M4}), $\phi$ is also a Markov trace.
The constants $\ta$ and $\bar\ta$ are now
\[
\ta=\pm\, {q^{-9/2}\over q^6+q^4+q^2+1},\qquad
\bar\ta=\pm\, {q^{21/2}\over q^6+q^4+q^2+1}
\]
that together with the minimal polynomial of $R$ in (\ref{R4})
\[
(R-Z)\,(R+q^6\,Z)\,(R-q^{10}\,Z)\,(R+q^{12}\,Z)=0,\qquad  Z=\pm\, q^{-9/2}
\]
give for the link invariant (\ref{ambientp}) the expression
\begin{eqnarray}
\al'(b_i^3)&=&
\pm\,t^{3/2}\,(1-t^3+t^5-t^6)\,
\al'(b_i^2)\nn\\
&&+t^6\,(1-t^2+t^3+t^5-t^6+t^8)\,
\al'(b_i)
\label{polpN4}\\
&&\mp\,t^{9/2}\,t^8\,(1-t+t^3-t^6)\,
\al'(1)-
t^{20}\,\al'(b_i^{-1})\nn
\end{eqnarray}
that coincides with (\ref{polN4}).

\mysection{Conclusion and generalization to $N$ arbitrary}
\label{conclusions}

The collection of these results obtained for $N=2,3,4$ indicates
that the equality
\[
\al'=\al_{\ts{\;\rm vertex}\; {\rm model}}\,,
\]
where $\al'$ is calculated with the link invariant recalled in the
introduction and $\al$ via vertex models as used in \cite{ADW} is also
true for arbitrary $N$. In the case of generic $N$ it is possible to
write the corresponding ambient isotopy link invariant $\al'$ with the
following formulae deduced by induction
\[
\ta=\pm\, {q^{{-(N-1)^2/2}}\over 
q^{2\,(N-1)}+q^{2\,(N-2)}+\cdots +1},\qquad
\bar\ta=\pm\, {q^{{(N-1)\,(N+3)/2}}\over 
q^{2\,(N-1)}+q^{2\,(N-2)}+\cdots +1},
\]
in addition to the $N$-th degree minimal polynomial of $R$ of
expression
\[
\prod^N_{i=1}\left(R-(-1)^{i+1}q^{{N\,(N-1)-(N-i+1)\,(N-i)}}\,Z\right),
\quad {\rm where}\quad  Z=\pm\, q^{{-(N-1)^2/2}}.
\]

\mysection{Relation between $<L>$ and other methods of obtaining 
$N$-state link invariants}
\label{six}

We discuss briefly at this point the connection between invariant
$<L>$ and other methods existing in the literature which are known to
provide Jones and Kauffman polynomials (i.e. the lowest $N$ link
invariants of the series) such as the Kauffman bracket approach \cite{Kau3}
\cite{Kau2}, or the work of Kirillov and Reshetikhin \cite{KR} which
reproduces the entire hierarchy of $N$-state link invariants.

\mysubsection{Kauffman bracket approach and state models} 

To discuss how the invariant $<L>$ is connected with 
Kauffman bracket approach we need to note first the following
fact: for any $R$-matrix with associated matrices $M_u$ and $M_d$ it
is possible to always define a Temperley-Lieb (TL) algebra, or
viewed in an equivalent manner, that any regular link invariant
$<L>$ constructed out of conditions (\ref{m})-(\ref{twist2}) comes
naturally equipped with a Temperley-Lieb algebra. The algebra is
formulated as follows: let $e=(e^a{}_c{}^b{}_d)$ with $a,b,c,d$ in the
index set $I$ be the matrix whose matrix elements are defined by
$e^a{}_c{}^b{}_d=M^{a\,b}\, M_{c\,d}$; then $e$ provides a
Temperley-Lieb algebra with generators $e_1,\ldots, e_{n-1}$ given by
\[
e_i=\un\oti\cdots\oti\underbrace{e}_{i,\, i+1}\oti\cdots\un
\]
and relations  
\begin{equation}
e_i^2=k\,e_i,\qquad e_i\,e_j=e_j\,e_i,\quad |i-j|\ge 2, \qquad 
e_i\,e_{i\pm 1}\,e_i=e_i
\label{TL}
\end{equation}
where the constant $k={{}\atop\epsfbox{zero.eps}}= M_{a\,b}\,M^{a\,b}$
and the indices $1\,\le i,\, j\,\le n-1$ are chosen so that all
relations above make sense. The proof of this statement is very
simple. The first relation comes from
$e^2={{}\atop\epsfbox{zero.eps}}\,e$, which obviously holds from the
definition of $e$ and condition (\ref{m}). The second relation is
clear when indices $i,\,j$ are sufficiently apart from each other. In
the case of the upper third expression, its lhs written in components
is $e^a{}_i{}^b{}_j\,e^j{}_q{}^c{}_f\,e^i{}_d{}^q{}_e= M^{a\,b}\,
M_{i\,j}\,M^{j\,c}\, M_{q\,f}\,M^{i\,q}\, M_{d\,e}=
M^{a\,b}\,M_{d\,e}\,\delta^c{}_f=e^a{}_d{}^b{}_e\,\delta^c{}_f$, which
equals the rhs when this is expressed in components too. The lower
relation is proved in a similar manner. (We note that there exist
another possible TL algebra defined not with $e$ but with matrix $f$
as $f^a{}_c{}^b{}_d=M^{b\,a}\, M_{d\,c}$ and that provides similar
relations to those in (\ref{TL}) but with $f$ instead of $e$ in the
definition of the generators $e_i$. However, the matrix representation
of this second TL algebras is related by a similarity transformation
to the one constructed with $e$ alone since $f=P\,e\,P$, where $P$ is
the permutation matrix, $P^a{}_c{}^b{}_d= \delta^a{}_d\,\delta^b{}_c$
(Kronecker deltas)).  

We use now this Temperley-Lieb algebra to formulate a state model
for the $N=2$ $R$-matrix given in (\ref{R2}) 
(from the two choices that we have for the constant $Z$ we take 
$Z=q^{-1/2}$ to work out the example) given by the expression
\[
R^a{}_c{}^b{}_d=A\,\delta^a{}_c\,\delta^b{}_d+ B\,M^{a\,b}\, M_{c\,d}, 
\]
were the indices are taken in the set $I=\{-1/2,\,1/2\}$, $M_u,\,M_d$
are the matrices written in (\ref{M2}) and $A,\,B$ are constants to be
determined. In this case, the minimal polynomial of $R$ (remember from
Section~\ref{four} that this is given by 
$q^{-1/2}\,R-q^{1/2}\,R^{-1}=q^{-1}-q$) fixes 
without ambiguity the constants $A,\,B$ and gives for $R$ above the 
formula $R=q^{-1/2}\,\un+q^{1/2}\,e$ where $e$ is, according to the
definition, the matrix
\begin{eqnarray*}
e=\pmatrix{
0&0&0&0\cr
0&-q&1&0\cr
0&1&-q^{-1}&0\cr
0&0&0&0\cr}.
\end{eqnarray*}
In pictorial form $R$ can be represented by the bracket
identity
\begin{equation}
{{}\atop\epsfbox{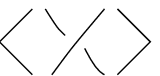}}=
q^{-1/2}\,{{}\atop\epsfbox{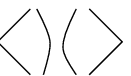}}
+q^{1/2}\,{{}\atop\epsfbox{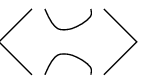}}
\label{bracket}
\end{equation}
that together with equations 
\begin{equation}
<{{}\epsfbox{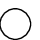}}\,L>=
-(q+q^{-1})\,<L>,\qquad
<{{}\epsfbox{zeroL.eps}}>=-(q+q^{-1})
\label{bracketcond}
\end{equation}
constitutes the bracket model for the matrix in (\ref{R2}). The first
equation in (\ref{bracketcond}) is obtained just applying
property $<{{}\epsfbox{zeroL.eps}}\,L>=M_{a\,b}\,M^{a\,b}\,<L>$
particularized for the case $N=2$.  Jones polynomial can now be
obtained from the bracket in the usual manner \cite{Kau2}: from
(\ref{bracket}) and (\ref{bracketcond}) derive the two following relations
\begin{eqnarray*}
&&<{{}\atop\epsfbox{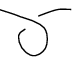}}>=-q^{-3/2}\,<{{}\atop\epsfbox{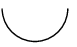}}>\\
&&<{{}\atop\epsfbox{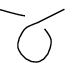}}>=-q^{3/2}\,<{{}\atop\epsfbox{line.eps}}>
\end{eqnarray*}
which indicate that the normalization of the bracket given by 
$V_L=(-q^{3/2})^{w(L)}<L>$
is an ambient isotopy invariant for oriented links. Here
$w(L)$, the writhe of the link $L$, is the regular isotopy invariant
defined by the equation
$w(L)=\sum_p \epsilon(p)$ where $p$ runs over all crossings of $L$ and
$\epsilon(p)$ is the sign of the crossing
$\epsilon\l({{}\atop\epsfbox{sigma1.eps}}\r)=1$,
$\epsilon\l({{}\atop\epsfbox{sigma-1.eps}}\r)=-1$. Also from 
(\ref{bracket}) derives the identity
\[
q^{-1/2}\,{{}\atop\epsfbox{stateR.eps}}
-q^{1/2}\,{{}\atop\epsfbox{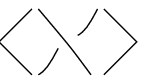}}=
(q^{-1}-q)\,{{}\atop\epsfbox{stateD.eps}}
\]
that written in terms of $V_L$ defined by the previous normalization
gives the Jones polynomial skein relation ($t=q^2$ as usual)
displayed in formula (\ref{polN2}), i.e.
\[
t^{-1}\,V{}_{{}_{{}\atop\epsfbox{sigma1.eps}}}-
t\,V{}_{{}_{{}\atop\epsfbox{sigma-1.eps}}}=
(t^{1/2}-t^{-1/2})\,V{}_{{}_{{}\atop\epsfbox{sigma0.eps}}}.
\]
This shows that matrices $M_u$ and $M_d$ are of relevance to associate
bracket identities to $R$-matrices thus making a direct link between
the invariant $<L>$ and the bracket approach introduced by
Kauffman. We have seen also with (\ref{bracket}) that they help to
formulate a state model for the matrix (\ref{R2}). The first fact, the
definition of bracket identities through $M_u$ and $M_d$, is common to the
whole series of $N$-state link invariants and therefore can be done
for generic $N$, but the definition of a state model when $N\ge 3$ is
a rather different subject. Let us see it by examining the case
$N=3$. With the expression of $M_u$ and $M_d$ given in (\ref{M3}) we
calculate the TL element $e$ which now results in the matrix
\begin{eqnarray}
e=\l(\begin{array}{ccccccccc}
0&0&0&0&0&0&0&0&0\\
0&0&0&0&0&0&0&0&0\\
0&0&q^2&0&-q&0&1&0&0\\
0&0&0&0&0&0&0&0&0\\
0&0&-q&0&1&0&-q^{-1}&0&0\\
0&0&0&0&0&0&0&0&0\\
0&0&1&0&-q^{-1}&0&q^{-2}&0&0\\
0&0&0&0&0&0&0&0&0\\
0&0&0&0&0&0&0&0&0
\end{array}\r)
\,\label{e3}
\end{eqnarray}
and that verifies $e^2=(q^2+1+q^{-2})\,e$. This element is useful to
write the following relation satisfied by matrix (\ref{R3}) with the
choice $Z=q^{-2}$
\[
R-R^{-1}=(q^{-2}-q^2)\,(\un-e),
\]
or equivalently
\begin{equation}
{{}\atop\epsfbox{stateR.eps}}-{{}\atop\epsfbox{stateRi.eps}}=
(q^{-2}-q^2)\,\l({{}\atop\epsfbox{stateD.eps}}-
{{}\atop\epsfbox{stateTL.eps}}\r).
\label{Dubrovnik1}
\end{equation}
This identity, together with 
\begin{eqnarray}
&&<{{}\epsfbox{zeroL.eps}}>=q^2+1+q^{-2}\nn\\
&&<{{}\atop\epsfbox{loop1.eps}}>=q^{-4}\,<{{}\atop\epsfbox{line.eps}}>
\label{Dubrovnik2}\\
&&<{{}\atop\epsfbox{loop2.eps}}>=q^{4}\,<{{}\atop\epsfbox{line.eps}}>\nn
\end{eqnarray}
obtained from (\ref{R3}) and (\ref{M3}) too define (a specialization
of) the Dubrovnik version of the Kauffman polynomial for unoriented
links. Yet again we see that matrices $M_u$ and $M_d$ are helpful to
formulate bracket identities but, unlike the $N=2$ case, the
$R$-matrix in (\ref{R3}) cannot be expanded in terms of the unit
matrix and $e$ only but requires the addition of extra terms. It is
clear then that there is an infinite set of state models that satisfy
the Dubrovnik version (\ref{Dubrovnik1}), (\ref{Dubrovnik2}) of the
Kauffman polynomial, but it is an open (and interesting) question to
decide whether there exists a state model which only involves
combination of Kronecker deltas, $M_{a\,b}$ and $M^{c\,d}$. If such
state model does exist the extra terms to add to matrices $\un$ and
$e$ should be of third order or bigger in the matrix elements of $M_u$
and $M_d$. This feature is common to all $N\ge 3$.

\mysubsection{Invariants of links from representations of $U_qsl(2)$} 

Kirillov and Reshetikhin have shown in ref. \cite{KR} that the link
invariants constructed through $N$-state vertex models do coincide with
the invariants derived from the spin $j$ representation of the quantum
enveloping algebra of $sl(2)$.  Since this paper reproduces the
$N$-hierarchy via $<L>$ we dedicate this section to explain the
connection between the invariant $<L>$ and the formalism of the
Russian authors.

Let us recall some properties of the quantum enveloping algebra
$U_qsl(2)$ that we need first.  $U_qsl(2)$ is the Hopf
algebra generated by elements $H,\,X^\pm$ with algebra relations
\begin{equation}
[H,X^\pm]=\pm 2\,X^\pm,\qquad [X^+,X^-]={q^H-q^{-H}\over q-q^{-1}}
\label{algebra}
\end{equation}
and coalgebra relations given by
\begin{eqnarray*}
&&\tri H=H\oti\un+\un\oti H,\qquad \tri X^\pm=q^{H/2}\oti X^\pm+
X^\pm\oti q^{-H/2}\\
&&\epsilon(H)=\epsilon(X^\pm)=0
\end{eqnarray*}
in addition to the antipode map $S(\cdot)$ with action
\begin{equation}
S(H)=-H,\qquad S(X^\pm)=-q^{\mp 1}\,X^\pm
\label{antipode}
\end{equation}
Here $q$ is considered as a generic real constant. This algebra
can be extended by adding an invertible element $w$ (which
obviously is not in $U_qsl(2)$) such that the element performs the
transformation \cite{KR}
\begin{equation}
w\,a\,w^{-1}=\ta\,S(a),\qquad{\rm for\; all}\; a\;{\rm in}\quad
U_q sl(2)
\label{weyl}
\end{equation}
where $\ta$ denotes now the linear anti-automorphism of action
\begin{equation}
\tau(H)=H,\qquad \tau(X^\pm)=X^\mp
\label{transposition} 
\end{equation}
on the generators of $U_q sl(2)$. A consequence of the existence of
$w$ is the `crossing symmetry' property exhibited by the universal
${\cal R}$-matrix of $U_q sl(2)$ (we do not review any property of this 
universal object here. Nevertheless, details about it can be found in
\cite{KR} as well)
\begin{equation}
{\cal R}=q^{-(H\oti H)/2}\sum_{n=0}^\infty {(1-q^2)^n\over [n]!}\,
q^{-n(n-1)/2}(q^{-H/2}X^+)^n\oti (q^{H/2}X^-)^n
\label{universal}
\end{equation}
(here $[x]=(q^x-q^{-x})/(q-q^{-1})$) which manifests in the formulae
\begin{eqnarray}
&&(\ta\oti {\rm id}) {\cal R}^{-1}=(w\oti {\rm id})\,{\cal R}\,
(w^{-1}\oti {\rm id})
\label{uni1}\\
&&({\rm id}\oti\tau){\cal R}=({\rm id}\oti w)\,{\cal R}^{-1}\,
({\rm id}\oti w^{-1})
\label{uni2}
\end{eqnarray}
These two formulae hold because the universal 
${\cal R}$-matrix of $U_q sl(2)$ is an
invertible operator which satisfies ${\cal R}^{-1}=
(S\oti {\rm id})\,{\cal R}$ and ${\cal R}=({\rm id}\oti S)\,{\cal R}^{-1}$.
Lowered to a representation $\pi^{j_1}\oti\pi^{j_2}$ of 
$U_q sl(2)^{\oti^2}$, equations (\ref{uni1}) and (\ref{uni2}) 
can be written as
\begin{eqnarray}
&&((R^{{j_1}{j_2}})^{-1})^{t_1}=(w^{j_1}\oti \un)\,R^{{j_1}{j_2}}
\,((w^{j_1})^{-1}\oti \un)
\label{cs1}\\
&&(R^{{j_1}{j_2}})^{t_2}=(\un\oti w^{j_2})\,(R^{{j_1}{j_2}})^{-1}
\,(\un\oti (w^{j_2})^{-1})
\label{cs2}
\end{eqnarray}
where we are using the notation $R^{{j_1}{j_2}}=\pi^{j_1}\oti\pi^{j_2}({\cal
R})$ and $w^j=\pi^j(w)$. Transpositions in the first and second space
in $V^{j_1}\oti V^{j_2}$ are indicated by $t_1$ and $t_2$. Equation
(\ref{cs1}) is precisely relation (2.10) in \cite{KR}.

We prove now that when $j_1=j_2=j$ the equalities (\ref{cs1}),
(\ref{cs2}) are exactly (\ref{twist1}), (\ref{twist2}) provided that
we make the identification $M_d=({w^j})^t$, i.e. that $M_d$ is the
transpose of the matrix associated to the element $w$ in the $j$
representation of $U_q sl(2)$. The proof of this statement is very
simple: relation (\ref{cs2}) is written in components as
\[
(R^{j\,j})^a{}_c{}^b{}_d=w^j{}{}^d{}_e\,{(R^{j\,j})^{-1}}^a{}_c{}^e{}_f\,
{(w^j)^{-1}}^f{}_b
\]
that expressed in terms of the $R$-matrix
$R^a{}_c{}^b{}_d=(R^{j\,j})^b{}_c{}^a{}_d$ looks as
\[
R^b{}_c{}^a{}_d=w^j{}{}^d{}_e\,{R^{-1}}^a{}_f{}^e{}_c\,
{(w^j)^{-1}}^f{}_b
\]
After reversing this relation it reads finally as
\[
{R^{-1}}^a{}_c{}^b{}_d={(w^j)^{-1}}^b{}_f\,
R^e{}_d{}^a{}_f\,
w^j{}^e{}_c
\]
Notice that this is precisely relation (\ref{twist2}) after the
substitution $M_d=({w^j})^t$ (remember the convention
adopted in this paper, namely $M_{a\,b}=(M_d)^a{}_b$ and 
$M^{a\,b}=(M_u)^a{}_b$). The proof that
(\ref{cs1}) is (\ref{twist1}) can be done in a similar fashion.

This much for the connection between $<L>$ and the formalism carried
out in \cite{KR} to construct link invariants from $U_q sl(2)$.
Let us see now as a practical example that indeed the
transpose of the matrix $\pi^j(w)$ when $j=1/2,1,3/2$ is the matrix
$M_d$ when $N=2,3,4$ as displayed in (\ref{M2}), (\ref{M3}) and  
(\ref{M4}), respectively (the correspondence is 
given by $j=(N-1)/2$ for generic $N$). {}From the algebra
relations (\ref{algebra}) and the Casimir operator of $U_q sl(2)$
given by
\[
C=\l({q^{(H+1)/2}-q^{-(H+1)/2}\over q-q^{-1}}\r)^2+X^-X^+=
\l({q^{(H-1)/2}-q^{-(H-1)/2}\over q-q^{-1}}\r)^2+X^+X^-
\]
is deduced the action of $H,\,X^\pm$ on the basis vectors
$\{e^j_m\}$, $m=-j,\ldots,j$ that span the $2\,j+1$ dimensional
representation $\pi^j$ of $U_q sl(2)$, i.e. 
\begin{equation}
\pi^j(H)\,e^j_m=2\,m\,e^j_m, \qquad 
\pi^j(X^\pm)\,e^j_m=([j\mp m]\,[j\pm m+1])^{1/2}\, e^j_{m\pm 1}
\label{j}
\end{equation}
This action is sufficient to calculate $\pi^j(w)$ up to a constant 
$\gamma_j$ that depends on the representation $j$ and that is of no
relevance for link invariance purposes. Indeed from (\ref{weyl}), 
(\ref{transposition}) and the antipode action (\ref{antipode})
it follows that 
\[
w\,H\,w^{-1}=-H,\qquad w\,X^\pm\,w^{-1}=-q^{\mp 1}\,X^\mp
\]
that together with the action (\ref{j}) allows to find the matrix
elements of $\pi^j(w)$ 
\[
w^j_{m\,m'}=<e^j_m|\pi^j(w)|\,e^j_{m'}>=(-1)^{j+m}\, q^{j+m}\,\gamma_j
\,\delta_{m,-m'}
\] 
In obtaining these matrix elements is important to notice that
$\pi^j(X^+)=(\pi^j(X^-))^t$, i.e. that we are truly under the
hypothesis expressed in the transposition (\ref{transposition}). Now
in the basis $\{e^j_{-j},\ldots,e^j_j\}$ and for $j=1/2,1,3/2$
we have that 
\begin{eqnarray*}
&&(\pi^{1/2}(w))^t=q^{1/2}\,\gamma_{1/2}\l(\begin{array}{cc}
0     &    -q^{1/2}      \\
q^{-1/2}  &    0
\end{array}\r),
\qquad 
(\pi^1(w))^t=q\,\gamma_{1}\l(\begin{array}{ccc}
0  & 0   &   q    \\
0  & -1  &    0   \\
q^{-1}  & 0   & 0
\end{array}\r)\\
&&(\pi^{3/2}(w))^t=q^{3/2}\,\gamma_{3/2}\l(\begin{array}{cccc}
0  & 0   & 0  &  q^{3/2}  \\
0  & 0   & -q^{1/2} &  0   \\
0  & q^{-1/2} & 0 & 0     \\
-q^{-3/2} & 0  &  0 & 0
\end{array}\r)
\end{eqnarray*}
result to be compared with matrices $M_d$ in (\ref{M2}), (\ref{M3}) and
(\ref{M4}). They do differ in an extra factor $q^j\,\gamma_j$ but we
remind again that $M_d$ can be determined up to an arbitrary
multiplicative constant and that the invariant $<L>$ is completely
independent of which the value of this constant is. This means
that the extra factor can be incorporated into $M_d$ and then
it is correct to identify
$M_d=(\pi^j(w))^t$ as we wanted to prove.

We mention that matrices (\ref{R2}), (\ref{R3}) and (\ref{R4}) are the
intertwining operator of the representation $V^j\oti V^j$ of
$U_qsl(2)$ in the cases $j=1/2$, $j=1$ and $j=3/2$, respectively. More
precisely they are $R=P^{j\,j}\,R^{j\,j}$, $P^{j\,j}$ the permutation
operator in $V^j\oti V^j$ and $R^{j\,j}=\pi^{j}\oti\pi^{j}({\cal R})$,
${\cal R}$ as in (\ref{universal}).

We conclude this section with a remark: the element $w$ defined by
relations (\ref{weyl}) and (\ref{transposition}) does exist for
$U_qsl(2)$ as we know. It also exists for other (not all) quantized
enveloping algebras after suitable modifications of
(\ref{transposition}) to allow transposition on the various generators
$H_i,\,X_i^\pm$. It follows then that the link
invariant construction introduced in \cite{KR} and its generalization
applies for (some) quantized enveloping algebras.  The
case of invariant $<L>$ is rather different. As explained in
Section~\ref{intro}, its construction requires merely to start with an
$R$-matrix (whether there exist matrices $M_u$ and $M_d$ for a given
$R$-matrix is a different question) and this matrix does not need to
come from quantized enveloping algebras necessarily.

\mysection{Final remarks}
\label{seven}

\begin{enumerate}
\item For the link invariant $<L>$ to exist it is not necessary that
the $R$-matrix entries satisfy the {\em charge conservation}
condition, i.e. that $R^a{}_c{}^b{}_d=0$ unless $a+b=c+d$. In the case
of the $R$-matrices (\ref{R2}), (\ref{R3}) and (\ref{R4}) considered
here this is the case since they all have this property, but the
existence of a non-trivial link invariant $<L>$ does not require such
condition.

\item Every $N$-state vertex model affords a $R$-matrix when the limit
$\lim_{u\to\infty} R(u)/\rh(u)$ is well-defined. This limit is
well-defined for all values of $\mu$ in the interval
$-1/2\,\le\,\mu\,\le\,1/2$. Along the paper we have worked the case
$\mu=1/2$, let us discuss now the link invariant $\al'$ when
$\mu=-1/2$. In this case and for all $N$ the corresponding $R$-matrix
is $P\,R\,P$, with $R$ the $R$-matrix corresponding to $\mu=1/2$ and
$P$ again the permutation matrix, $P^a{}_c{}^b{}_d=
\delta^a{}_d\,\delta^b{}_c$. If ($R$, $M_u$) is a solution of eqs
(\ref{m})-(\ref{twist2}) with associated link invariant $<L>$, then
($P\,R\,P$, $M_u^t$) is also a solution with associated link invariant
$<L'>$, where $L'$ is as $L$ but with the braid strands closed on the
opposite side of the plane.  Now $<L>=<L'>$ and consequently
$\al'|_{\mu=1/2}=\al'|_{\mu=-1/2}$, i.e. the link invariants for
$\mu=1/2, -1/2$ are the same.  Regarding the case in which
$-1/2\,<\,\mu\,<\,1/2$, there exist matrices $M_u$ and $M_d$ for each
$N$ if and only if $q$ is the following root of unity
$q^{2(N-1)}=1$. For these roots the minimal polynomial of $R$
collapses to $R^2=Z^2$ ($Z$ the values listed in the conclusions
Section~\ref{conclusions} restricted to the mentioned roots of unity)
and the associated link invariant skein relation is
$\al'\l({{}\atop\epsfbox{sigma1.eps}}\r)=
\al'\l({{}\atop\epsfbox{sigma-1.eps}}\r)$. This is an invariant
related to the number of components of the link and depends on $N$
simply because $\al'{{}\atop\epsfbox{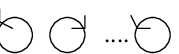}}$ of an arbitrary
number of unknots depends on $N$.
\end{enumerate}
\begin{figure}
\[{{}\atop\epsfbox{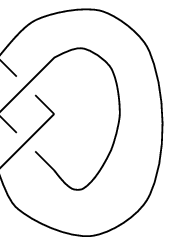}}\qquad\qquad {{}\atop\epsfbox{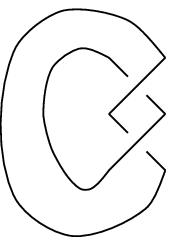}}\]
\caption{An example of $L$ and $L'$} 
\end{figure}


%
\def\section{\subsection}

\end{document}